\documentclass[12pt]{article}
\usepackage{amssymb,amsmath,epsfig,amscd,eucal}

\newtheorem{thm}{Theorem}[section]
\newtheorem{defn}[thm]{Definition}
\newtheorem{lem}[thm]{Lemma}
\newtheorem{prop}[thm]{Proposition}

\newtheorem{rem}[thm]{Remark}
\newtheorem{cor}[thm]{Corollary}

\newenvironment{pf}{\par\medskip\noindent{\em Proof. }}{\hfill $\square$\par\medskip}
\newenvironment{pfof}[1]{\par\medskip\noindent{\em Proof of #1. }}{\hfill $\square$\par\medskip}
\newenvironment{altpfof}[1]{\par\medskip\noindent{\em Alternative proof of #1. }}{\hfill $\square$\par\medskip}

\newcommand{\R}{\mathbb{R}}
\newcommand{\Z}{\mathbb{Z}}

\title{Elementarily free groups are subgroup separable}
\author{Henry Wilton}
\date{16th November 2005}

\begin{document}
\maketitle

\begin{abstract}
Elementarily free groups are the finitely generated groups with the
same elementary theory as free groups.  We prove that elementarily
free groups are subgroup separable, answering a question of Zlil
Sela.
\end{abstract}

\emph{Limit groups} arise naturally in the study of the set of
homomorphisms to free groups and, in the guise of \emph{fully
residually free groups}, have long been studied in connection with
the first-order logic of groups: see for, example, \cite{KM98a},
\cite{KM98b} and \cite{Rem89}. Indeed, limit groups turn out to be
precisely the groups with the same existential theory\footnote{The
existential theory of $G$ is the set of sentences in the elementary
theory that only use the existential quantifier $\exists$.} as a
free group \cite{Rem89}. The name limit group was introduced by Zlil
Sela in his solution to the Tarksi Problem (see \cite{Se1},
\cite{Se2} \emph{et seq.}), wherein he characterized the finitely
generated groups with the same elementary theory\footnote{The
elementary theory of a group $G$ is the set of first-order sentences
that are true in $G$.} as a free group. Such group are called
\emph{elementarily free}. Sela asked if limit groups were subgroup
separable in \cite{SeQ}.

A group $G$ is \emph{subgroup separable} (or \emph{LERF}) if any
finitely generated subgroup $H$ and any element $g\notin H$ can be
distinguished in a finite quotient.  The generalized word problem
is soluble for subgroup separable groups.  Historically, it has
been of interest in topology because, in certain circumstances, an
immersion into a space with subgroup separable fundamental group
can be lifted to an embedding in a finite cover (see, for example,
\cite{Th}).

The first significant result concerning subgroup separability is
due to M.~Hall in \cite{Hall49}, who demonstrated the property for
free groups.  This was generalized by R.~G.~Burns \cite{Bu} and
N.~S.~Romanovskii \cite{Ro}, who showed that a free product of
subgroup separable groups is subgroup separable.

These results were all proved using algebraic methods, but a more
topological approach was developed by J.~Hempel in \cite{He76},
J.~R.~Stallings in \cite{St} and P.~Scott in \cite{Sc1} (see also
\cite{Sc2}); Scott used hyperbolic geometry to prove that surface
groups are subgroup separable.

A.~M.~Brunner, R.~G.~Burns and D.~Solitar \cite{BBS} extended
Scott's theorem by showing that any amalgamated product of free
groups along a cyclic subgroup is subgroup separable.  The
topological thread was continued by, among others, M.~Tretkoff \cite
{Tr} and G.~Niblo \cite{Ni}. More recently, D.~Long and A.~Reid
\cite{LR} adapted Scott's approach to show that geometrically finite
subgroups of certain hyperbolic Coxeter groups are separable.

In preparation for our main theorem, we give a purely topological
proof of Scott's theorem (see corollary \ref{Surface groups}), as
well as proofs of the results of Hall, Burns and Romanovskii.

\begin{thm}[Theorem \ref{Main theorem}]
Elementarily free groups are subgroup separable.
\end{thm}

The proof of theorem \ref{Main theorem} is most closely prefigured
by the work of Rita Gitik and Daniel Wise.  Wise \cite{Wi}
classified the graphs of free groups with cyclic edge groups that
are subgroup separable. Gitik showed in \cite{Gi1} that if $G$ is
subgroup separable, $F$ is free and $f\in F$ has no proper roots
then $G*_{\langle f\rangle}F$ is also subgroup separable. Section
\ref{Less direct} explains how one can combine this result with
recent work of Martin Bridson, Michael Tweedale and myself on
positive-genus towers to provide a second, less direct proof of
theorem \ref{Main theorem}.

This paper is organized as follows. In section 1 we introduce the
definition of subgroup separability and some basic properties and
reformulations. In section 2 we develop our basic tools, the notions
of pre-cover and elevation. The fundamental Stallings' Principle of
section 3 gives a criterion for completing a pre-cover to a cover,
which we use to give topological proofs of the theorems of Hall,
Burns and Romanovskii.  We anticipate that this formulation of
Stallings' Principle will prove useful in other contexts.  We show
how assumptions on homology make it possible to apply Stallings'
Principle and provide a topological proof of Scott's result that
surface groups are subgroup separable.  A detailed argument in
section 4 shows how to replace infinite pre-covers by finite ones.
In section 5 we introduce elementarily free groups and prove that
they are subgroup separable.

The arguments of the current article are framed with an attack on
the general problem of whether all limit groups are subgroup
separable in mind.  These techniques are, however, ill-suited to the
natural, stronger question of whether all finitely generated
subgroups of limit groups are closed in the pro-free topology.

\section{Subgroup separability}

\begin{defn}
Let $G$ be a group.  A subgroup $H\subset G$ is \emph{separable}
if it is an intersection of finite-index subgroups of $G$. Call
$G$ \emph{subgroup separable} if every finitely generated subgroup
is separable.
\end{defn}

Note that a subgroup $H$ is separable if and only if, for any
$g\notin H$, there exists a finite-index subgroup $K$ of $G$ so
that $H\subset K$ and $g\notin K$.  Here are some other useful
algebraic reformulations.

\begin{lem}\label{Algebraic rephrasings}
Let $G$ be a group and $H$ a subgroup.  The following are equivalent.
\begin{enumerate}
\item $H$ is separable.
\item For every $g\notin H$ there exists a finite-index subgroup $K\subset
G$
so that $g\notin HK$
\item For every $g\notin H$ there exists a homomorphism $f:G\to Q$
to a finite group so that $f(g)\notin f(H)$.
\end{enumerate}
\end{lem}
\begin{pf}
That 1 implies 2 is trivial.  Given $K$ as in 2, we can construct a canonical
normal subgroup
$$
N=\bigcap_{g\in G} K^g.
$$
It is of finite index in $G$ because if $g_1,\ldots,g_n$ are coset representatives
for $K$, $N=\cap_iK^{g_i}$.  Setting $Q=G/N$ shows that 2 implies 3.

To deduce 1 from 3, take $f^{-1}(f(H))$ to be the required finite-index subgroup.
\end{pf}

We will also make use of some equivalent topological notions. We
work in the category of combinatorial complexes and combinatorial
maps.

\begin{lem}\label{Topological rephrasings}
Let $(X,x)$ be a based, connected complex, $H\subset\pi_1(X,x)$ a
subgroup and $p:(X^H,x')\to (X,x)$ the covering corresponding to
$H$. The following are equivalent.
\begin{enumerate}
\item $H$ is separable.
\item For every finite subcomplex $\Delta\subset X^H$ there exists an intermediate, finite-sheeted
covering
$$
X^H\to\bar{X}\to X
$$
so that $\Delta$ embeds in $\bar{X}$.
\item For any $g\notin H$, there exists a finite-sheeted covering
$$
(\hat{X},\hat{x})\to(X,x)
$$
so that, for every $h\in H$, the end-point of the lift of any
(based) representative of $g$ to $\hat{X}$ doesn't coincide with
the end-point of the lift of any representative of $h$.
\end{enumerate}
\end{lem}
\begin{pf}
We start by showing how 2 follows from 1.  Let $\{\sigma^i_j\}$ be
the set of cells of $\Delta$, partitioned so that $\sigma^i_j$ and
$\sigma^{i'}_{j'}$ have the same image in $X$ if and only if
$i=i'$.  Let $x^i_j$ be the barycentre of $\sigma^i_j$.  For each
$i$ and $j$ fix a continuous path $\alpha_i$ from $x'$ to $x^i_j$
and, for each $i$ and distinct $j$ and $j'$, fix a continuous path
$\beta^i_{j,j'}$ from $x^i_j$ to $x^i_{j'}$.  The concatenation
$$
(p\circ\alpha^i_j)\cdot
(p\circ\beta^i_{j,j'})\cdot(p\circ\alpha^i_j)^{-1}
$$
defines an element $g^i_{j,j'}\in\pi_1(X)\smallsetminus H$.  Let $K^i_{j,{j'}}$
be a finite-index subgroup of $\pi_1(X)$ containing $H$ but not
$g^i_{j,j'}$.  The finite-sheeted covering $(\bar{X},\bar{x})\to
(X,x)$ with
$$
\pi_1(\bar{X})=\bigcap_{i,j,j'} K^i_{j,j'}
$$
is as required.  For suppose $\Delta$ does not embed.  Then the
images of two cells, $\sigma^i_j$ and $\sigma^{i'}_{j'}$, and
hence their barycentres, coincide in $\bar{X}$.  If this happens
then their images coincide in $X$, so $i=i'$ and
$g^i_{j,j'}\in\pi_1(\bar{X})$, a contradiction.

To see that 2 implies 3, fix a representative loop $\gamma$ for
$g$ (it doesn't matter which) and take $\Delta$ to be the image of the lift
of $\gamma$ to $X^H$. Then the lift of $\gamma$ to $\bar{X}$ is not closed,
whereas the lift of any representative of $h\in H$ to $\bar{X}$ is
closed.

Given $\hat{X}$ as in 3, set $K=\pi_1(\hat{X})\subset\pi_1(X)$.
Then if $g\in hK$ for $h\in H$ it follows that $gh^{-1}\in K$, so
the end-point of any lift of a representative of $g$ coincides
with the end-point of any lift of a representative of $h$, a
contradiction.  So 3 implies 1.
\end{pf}

It is immediate that the property of being subgroup separable is
closed under taking subgroups.  The following lemma will also be
useful.

\begin{lem}\label{Elementary properties}
Consider a subgroup separable group $G$.  If $G$ is a finite-index
subgroup of a group $G'$, then $G'$ is subgroup separable.
\end{lem}
\begin{pf}
Replacing $G$ with the intersection of its conjugates, it can be
assumed to be normal in $G'$.   Let $H\subset G'$ be a finitely
generated subgroup and $\gamma\in G'\smallsetminus H$.   If
$\gamma\notin HG$ then the result is immediate, so assume
$\gamma=hg$ for $h\in H$ and $g\in G\smallsetminus G\cap H$.  Since
$G$ is subgroup separable and $G\cap H$ is finitely generated, there
exists a finite-index normal subgroup $K\subset G$ so that $g\notin
(G\cap H)K=G\cap (HK)$. Therefore $\gamma\notin HK$, as required.
\end{pf}

\section{Elevations and pre-coverings}

\subsection{Graphs of spaces}

A \emph{graph of spaces} $\Gamma$ consists of:
\begin{enumerate}
\item a set $V(\Gamma)$ of connected spaces, called \emph{vertex spaces};
\item a set $E(\Gamma)$ of connected spaces, called \emph{edge spaces};
\item for each edge space $e\in E(\Gamma)$ a pair of $\pi_1$-injective continuous \emph{edge maps}
$$
\partial_\pm^e:e\to\bigsqcup_{V\in V(\Gamma)}V.
$$
When the edge in question is unambiguous, we often suppress the superscript and refer to $\partial_\pm^e$ simply as $\partial_\pm$.
\end{enumerate}

The associated topological space $|\Gamma|$ is defined as the
quotient of
$$
\bigsqcup_{V\in V(\Gamma)} V \sqcup \bigsqcup_{e\in E(\Gamma)}(e\times [-1,+1])
$$
obtained by identifying $(x,\pm 1)$ with $\partial_\pm^e(x)$ for
each edge space $e$ and every $x\in e$.  We will usually assume that $|\Gamma|$
is connected.  If $X=|\Gamma|$ we will
often say that \emph{$\Gamma$ is a graph-of-spaces decomposition
for $X$}, of just that \emph{$X$ is a graph of spaces}.  The
\emph{underlying graph} of $\Gamma$ is the abstract graph given by
replacing every vertex and edge space of $\Gamma$ by a point.  If
$X$ is the topological space associated to the graph of spaces
$\Gamma$, the underlying graph is denoted $\Gamma(X)$. Note that
there is a natural $\pi_1$-surjective map $\phi:X\to\Gamma(X)$.

Given a graph of spaces $X$, consider a subgraph
$\Gamma'\subset\Gamma(X)$.  The corresponding graph of spaces
$X'=\phi^{-1}(X)$ has a natural inclusion $X'\hookrightarrow X$ and
is called a \emph{sub-graph of spaces} of $X$.  If $X$ and $Y$ are
graphs of spaces, a map $f:X\to Y$ is a \emph{map of graphs of spaces} if, whenever $Y'$ is a sub-graph of spaces of $Y$, the pre-image $f^{-1}(Y')$ is a sub-graph of spaces of $X$.

The fundamental group of a graph of spaces is naturally a graph of
groups by the Seifert--van Kampen Theorem.  For more on graphs of
spaces, graphs of groups and Bass--Serre theory see \cite{SW} and
\cite{S77}.

\subsection{Elevations to covers}

Elevations are a natural generalization of lifts, and were introduced by Wise in \cite{Wi}.

\begin{defn}
Consider a continuous map of connected based spaces $f:(A,a)\to
(B,b)$ and a covering $B'\to B$.  An \emph{elevation} of $f$ to
$B'$ consists of a connected covering $p:(A',a')\to (A,a)$ and a
lift $f':A'\to B'$ of  $f\circ p$ so that for every intermediate
covering
$$
(A',a')\to(\bar{A},\bar{a})\stackrel{q}{\to}(A,a)
$$
there is no lift $\bar{f}$ of $f\circ q$ to $B'$ with
$\bar{f}(\bar{a})=f'(a')$.

Elevations $f'_1:A'_1\to B'$ and $f'_2:A'_2\to B'$ are
\emph{isomorphic} if there exists a homeomorphism $\iota:A'_1\to
A'_2$, covering the identity map on $A$, such that
$$
f'_1=f'_2\circ\iota.
$$
\end{defn}

In practice, we will often abuse notation and refer to just the
lift $f'$ as an elevation of $f$.  The next lemma follows from
standard covering-space theory.  See, for example, proposition
1.33 of \cite{Ha02}.

\begin{lem}\label{Covering space theory}
Fix a lift $b'\in B'$ of $b=f(a)$.  Consider the covering $p:(A',a')\to(A,a)$ such that
$$
\pi_1(A',a')=f_*^{-1}(\pi_1(B',b')).
$$
The composition $f\circ p$ admits a lift $f':A'\to B'$ to $B'$, which is an elevation of $f$.
\end{lem}

The \emph{degree} of the elevation $f':A'\to  B'$, denoted
$\deg(A')$, is the \emph{conjugacy class} of the subgroup
$\pi_1(A')\subset\pi_1(A)$.  If $\pi_1(A')$ is of finite index in
$\pi_1(A)$ then $f'$ is called \emph{finite-degree}; otherwise, $f'$
is \emph{infinite-degree}.

\begin{rem}
Consider a map $f:X\to Y$, a covering $Y'\to Y$ and an elevation $f':X'\to
Y'$ of $f$.  Let
$$
Y'\to\bar{Y}\to Y
$$
be an intermediate covering.  Then there exists a unique elevation
$\bar{f}:\bar{X}\to\bar{Y}$ of $f$ such that $X'\to X$ factors
through $\bar{X}\to X$ and

\[ \begin{CD}
{X'} @>f'>>  {Y'} \\
@VVV    @VVV \\
 {\bar{X}} @>\bar{f}>>  {\bar{Y}}\\
\end{CD}\]
commutes, determined by the requirement that
$\pi_1(\bar{X})=f_*^{-1}(\pi_1(\bar{Y}))$. We say $f'$
\emph{descends to $\bar{f}$}.
\end{rem}

\begin{rem}
If $X$ has a graph-of-spaces decomposition $\Gamma$ and $X'\to X$
is a covering space then $X'$ inherits a graph-of-spaces
decomposition $\Gamma'$, with vertex spaces the connected
components of the pre-images of the vertex spaces of $X$ and edge
spaces and maps given by all the elevations of the edge maps to
the vertex spaces of $X'$, up to isomorphism.
\end{rem}

\subsection{Elevations of embeddings}

These two lemmas, concerning elevations of embeddings, will prove
useful later.

\begin{lem}\label{Disjoint images for elevations of embeddings}
Let $f:X\to Y$ be an embedding and $Y'\to Y$ a connected covering.
If $f'_1,f'_2$ are non-isomorphic elevations of $f$ to $Y'$ then
their images are disjoint.
\end{lem}
\begin{pf}
Suppose there exists $y'\in\mathrm{im} f'_1\cap\mathrm{im} f'_2$.
Let $y\in Y$ be the image of $y'$ and set $x=f^{-1}(y)\in X$. By
lemma \ref{Covering space theory}, both $f'_1$ and $f'_2$ have the
same domain, namely the covering space $(X',x')\to (X,x)$
corresponding to $f_*^{-1}\pi_1(Y',y')$.
 The identity map from $X'$ to itself realizes an isomorphism between
 $f'_1$ and $f'_2$.
\end{pf}

\begin{lem}\label{Elevations of injective maps are injective}
If $f:X\to Y$ is an embedding and $f':X'\to Y'$ is an elevation then
$f'$ is an embedding.
\end{lem}
\begin{pf}
Suppose $x'_1,x'_2\in X'$ and $f'(x'_1)=f'(x'_2)=y'$.  Let $y\in Y$
be the image of $y'$ and set $x=f^{-1}(y)$.  Fix a path
$\gamma:[0,1]\to X'$ from $x'_1$ to $x'_2$.  Then $f'\circ\gamma$
represents an element of $\pi_1(Y',y')$.
 By lemma \ref{Covering space theory} we have that $f_*^{-1}(f'\circ\gamma)\in\pi_1(X',x'_i)$
 for either $i$, so in fact $\gamma$ is a loop and $x'_1=x'_2$.
\end{pf}

\subsection{Pre-covers}

Gitik uses a notion of pre-cover extensively in \cite{Gi1}.  Our
pre-covers are analogous to hers.  Pre-coverings are the
graph-of-spaces versions of Stallings' graph immersions in
\cite{St}.

\begin{defn}
Let $X$ and $\bar{X}$ be graphs of spaces ($\bar{X}$ is not
assumed connected). A \emph{pre-covering} is a locally injective
map $\bar{X}\to X$ that maps vertex spaces and edge spaces of
$\bar{X}$ to vertex spaces and edge spaces of $X$ respectively,
and restricts to a covering on each vertex space and each edge
space.  Furthermore, for each edge space $\bar{e}$ of $\bar{X}$
mapping to an edge space $e$ of $X$, the diagram of edge maps
\[ \begin{CD}
{\bar{e}} @>{\bar{\partial}_\pm}>>  {\bar{V}_\pm} \\
@VVV    @VVV \\
 {e} @>{\partial_\pm}>>  {V_\pm}\\
\end{CD}\]
is required to commute.  The domain $\bar{X}$ is a
\emph{pre-cover}.

The pre-covering $\bar{X}\to X$ is \emph{finite-sheeted} if the
pre-image of every point of $X$ is finite.
\end{defn}

\begin{rem}
All the edge maps of $\bar{X}$ are elevations of edge maps of $X$
to the vertex spaces of $\bar{X}$.  An elevation of an edge map of
$X$ to a vertex of $\bar{X}$ that isn't an edge map of $\bar{X}$
is called \emph{hanging}.  If none of the elevations are hanging
then $\bar{X}$ is in fact a cover.
\end{rem}

We will be interested in ways of completing pre-coverings to
genuine coverings in the spirit of Stallings \cite{St}.  First, we
find a canonical way of doing so that doesn't preserve finiteness.

\begin{prop}\label{Canonical completion of pre-covers}
Let $X$ be a graph of spaces and $\bar{X}\to X$ a pre-covering with
$\bar{X}$ connected. Then $\pi_1(\bar{X})$ injects into $\pi_1(X)$
and, furthermore , there exists a unique embedding
$\bar{X}\hookrightarrow\tilde{X}$ into a connected covering
$\tilde{X}\to X$ such that $\pi_1(\tilde{X})=\pi_1(\bar{X})$ and
$\tilde{X}\to X$ extends $\bar{X}\to X$.  This is called the
\emph{canonical completion} of $\bar{X}\to X$.
\end{prop}

Disconnected pre-covers also have canonical completions, given by
the disjoint union of the canonical completions of their connected components.

To prove proposition \ref{Canonical completion of pre-covers} we
will adapt some ideas from Bass--Serre theory. If $X$ is a graph of
spaces, a path $\gamma:I\to X$ is \emph{reduced} if it can't be
homotoped (relative to its end-points) off any vertex or edge space
of $X$ that it intersects. It is clear that every path is homotopic
to a reduced path.

\begin{lem}\label{Images of reduced paths are reduced}
If $\bar{X}\to X$ is a pre-covering and $\bar{\gamma}:I\to
\bar{X}$ is reduced then the image $\gamma$ of $\bar{\gamma}$ in
$X$ is also reduced.
\end{lem}
\begin{pf}
Suppose $V$ is a vertex of $X$ that $\gamma$ can be homotoped off.
Then lifting the homotopy to the pre-images of $V$ in $\bar{X}$,
it follows that $\bar{\gamma}$ can be homotoped off all the
pre-images of $V$.  So $\bar{\gamma}$ isn't reduced.
\end{pf}

We can now prove proposition \ref{Canonical completion of
pre-covers}.

\begin{pfof}{proposition \ref{Canonical completion of pre-covers}}  Let $p:\bar{X}\to X$ be a connected pre-covering.  First we show that $p_*:\pi_1(\bar{X})\to\pi_1(X)$
is an injection.  Let $\bar{\gamma}$ be a reduced loop in $\bar{X}$ and suppose that the image $\gamma$ of $\bar{\gamma}$ in $X$ is null-homotopic.  Since
$\gamma$ is reduced by lemma \ref{Images of reduced paths are reduced}, $\gamma$ is contained in a vertex space $V$ of $X$ and so $\bar{\gamma}$ is contained in a vertex space $\bar{V}$ of $\bar{X}$.  But $\bar{V}\to V$ is a covering map, so $\bar{X}$ is null-homotopic.

Let $\tilde{X}\to X$ be the covering corresponding to $p_*\pi_1(\bar{X})\subset\pi_1(X)$.  Then, by standard covering space theory, the map $p$ lifts to a map $\tilde{p}:\bar{X}\to\tilde{X}$.  It remains to show that $\tilde{p}$ is injective.

Suppose $x,y\in\bar{X}$ with $\tilde{p}(x)=\tilde{p}(y)$, and let
$\bar{\gamma}:I\to\bar{X}$ be a reduced path with end-points $x$
and $y$.  Then $\gamma=p\circ\bar{\gamma}$ is a loop in $X$. Since
$\tilde{p}$ is $\pi_1$-surjective, it can be assumed that $\gamma$
is null-homotopic.  But $\gamma$ is reduced, so in fact $\gamma$
is contained in a single vertex space $V$ of $X$.  Hence,
$\bar{\gamma}$ is contained in a single vertex space $\bar{V}$ of
$\bar{X}$.  Let $\tilde{V}=\tilde{p}(\bar{V})$, which is a
vertex space of $\tilde{X}$.  Then
$\tilde{\gamma}=\tilde{p}\circ\bar{\gamma}$ is a null-homotopic
loop in $\tilde{V}$, and the homotopy to a point lifts to
$\bar{V}$.  So $x=y$.
\end{pfof}

\subsection{Elevations to pre-covers}

Proposition \ref{Canonical completion of pre-covers} makes it possible to define elevations to pre-covers.

\begin{defn}\label{Elevations to pre-covers}
Let $f:X\to Y$ be a map of graphs of spaces and $\bar{Y}\to Y$ a
pre-covering.  Let $\tilde{Y}\to Y$ be the canonical completion of
$\bar{Y}\to Y$. Consider any finite, pairwise non-isomorphic
collection of elevations $\tilde{f}_i:\tilde{X}_i\to\tilde{Y}$ of
$f$. For each $i$, let $\bar{X}_i=\tilde{f}^{-1}(\bar{Y})\subset
\tilde{X}_i$.  Note that $\bar{X}=\sqcup_i\bar{X}_i$ is a pre-cover
of $X$ since $f$ is a map of graphs of spaces. Suppose that the
pre-covering $\bar{X}\to X$ is non-empty and can be extended to a
connected covering of $X$.  Then the coproduct
$$
\bar{f}=\sqcup_i\tilde{f}_i|_{\bar{X}_i}:\bar{X}\to\bar{Y}
$$
is called an \emph{elevation of $f$ to $\bar{Y}$}.

When $\bar{X}\to X$ is a genuine covering of $X$ the elevation $\bar{f}$ is called \emph{full}.
\end{defn}

In most examples of elevations to pre-covers in this article,
$X\cong S^1$ and $\bar{X}$ consists of a finite collection of closed
intervals.  In this case, $\bar{X}\to X$ can be realized as a
restriction of the universal covering $\R\to S^1$.

Given some maps of graphs of spaces $f_i:X\to Y$ and a pre-covering $\bar{Y}\to Y$ we will often be interested in a finite set of elevations $\{\bar{f}_j:\bar{X}_j\to
\bar{Y}\}$, where each $\bar{f}_j$ is an elevation of some $f_i$.
Such a set is called \emph{disjoint} if the images of distinct
$\bar{f}_j$ and $\bar{f}_k$ are disjoint.

\begin{defn}
Consider a map $f:S^1\to X$ and a pre-covering $\bar{X}\to X$ with canonical
completion $\tilde{X}$.  Let $\bar{f}:\bar{S}^1\to\bar{X}$
be an elevation that arises as a restriction of the coproduct of elevations
$$
\tilde{f}=\sqcup_i\tilde{f}_i:\tilde{S}^1=\sqcup_i\tilde{S}^1_i\to\tilde{X}
$$
as in definition \ref{Elevations to pre-covers}.  Consider $x\in
\partial\bar{S}^1$ and let $C\subset\tilde{S}^1$ be the closure of
the component of the complement of $\bar{S}^1$ such that
$x\in\partial C$.  The interval $C$ can be identified with a
subinterval of $[0,\infty)$ so that $x\equiv 0$. For some unique
edge space $\tilde{e}$ of $\tilde{X}$ and sufficiently small
$\epsilon>0$, $\tilde{f}([0,\epsilon])\subset
\tilde{e}\times[-1,+1]$.  We say that \emph{$\bar{f}$ extends to
$\tilde{e}$ at $x$}.  The elevation $\bar{f}:\bar{S}^1\to\bar{X}$ is
\emph{diverse} if $\bar{f}$ extends to distinct edge spaces of the
canonical completion at distinct points $x,y\in\partial\bar{S}^1$.

More generally, let $\{f_i:c_i\to X\}$ be a collection of maps from
circles $c_i$ to $X$ and consider a set of elevations
$\{\bar{f}_j:\bar{c}_j\to\bar{X}\}$ of the $f_i$.  The set
$\{\bar{f}_j\}$ is \emph{diverse} if $\bar{f}_j$ and $\bar{f}_k$
respectively extend to distinct edge spaces at distinct
$x\in\partial\bar{c}_j$ and $y\in\partial\bar{c}_k$.
\end{defn}

We will want to extend the pre-covering $\bar{X}\to X$ so that the
elevation $\bar{f}:\bar{S}^1\to\bar{X}$ extends to a full elevation.
The next proposition gives one example when this can be done.

\begin{prop}\label{Subtly extending graphs of spaces}
Suppose $\bar{X}\to X$ is a pre-covering such that all hanging
elevations of edge maps have simply connected domains. Let
$f:S^1\to X$ be a combinatorial map and $\bar{f}:\bar{S}^1\to
\bar{X}$ a diverse elevation of $f$ to $\bar{X}$. Let
$$
\hat{S}^1\to S^1
$$
be a finite-sheeted covering of $S^1$ and suppose there exists an embedding $\bar{S}^1\hookrightarrow\hat{S}^1$ so that $\hat{S}^1\to S^1$ extends $\bar{S}^1\to S^1$. Then there exists a pre-covering $\hat{X}\to X$ extending $\bar{X}\to X$ so that $\bar{f}:\bar{S}^1\to\bar{X}$ extends to a full elevation
$$
\hat{f}:\hat{S}^1\to\hat{X}
$$
of $f$.
\end{prop}
\begin{pf}
Let $C$ be the closure of a component of $\hat{S}^1\smallsetminus\bar{S}^1$.
Without loss identify $C\equiv [-1,1]$.  Consider the canonical
completion $\tilde{X}\to X$ of $\bar{X}\to X$.  There is a unique
lift $\tilde{f}_+:C\to\tilde{X}$ of $f$ to $\tilde{X}$ so that
$\tilde{f}_+(1)=\bar{f}(1)$.  Let $\epsilon>0$ be maximal such that $\tilde{f}_+(\epsilon-1)$
lies in a vertex space.    Let $X'\subset\tilde{X}$
be the pre-cover consisting of $\bar{X}$ together with the vertex spaces and edge spaces containing
$\tilde{f}_+((\epsilon-1,1))$.  Let $\tilde{e}_+$ be the edge space of
$\tilde{X}$ so that $\tilde{f}_+([-1,\epsilon-1])\subset \tilde{e}_+\times[-1,1]$
and without loss assume that $\tilde{f}_+(\epsilon-1)\in\tilde{e}_+\times\{+1\}$.

Similarly, there is a unique lift $\tilde{f}_-:C\to\tilde{X}$ of $f$ so that
$\tilde{f}_-(-1)=\bar{f}(-1)$.  Then $\tilde{f}_-([-1,\epsilon-1])\subset
\tilde{e}_-\times[-1,1]$ for some unique edge space $\tilde{e}_-$ of $\tilde{X}$.
Without loss, assume $\tilde{f}_-(-1)\in\tilde{e}_-\times\{-1\}$.

The edge spaces $\tilde{e}_+$ and $\tilde{e}_-$ are both simply connected
covers of some edge space $e$ of $X$.  There exists some unique covering
transformation $\tau:\tilde{e}_+\to\tilde{e}_-$ such that, whenever $x\in(\epsilon-1,-1)$
with $\tilde{f}_+(x)\in\tilde{e}_+\times\{\frac{1}{2}\}$,
$$
\tau\circ\tilde{f}_+(x)=\tilde{f}_-(x).
$$
Now let $\hat{X}$ be the pre-cover given by $X'$ together with the
additional edge space $\tilde{e}_+$; the additional edge maps are
$\partial^{\tilde{e}_+}_+$ and $\partial^{\tilde{e}_-}_-\circ\tau$.

Since $\bar{f}$ is diverse, this can be done for every of
component of $\hat{S}^1\smallsetminus \bar{S}^1$.
\end{pf}

\begin{rem}\label{Diversity is useful}
Consider a finite, diverse set of elevations $\{\bar{f}_j\}$ to $\bar{X}$
of various maps $f_i:S^1\to X$, so that each $\bar{f}_j$ satisfies
the conditions of proposition \ref{Subtly extending graphs of
spaces}.  Then proposition \ref{Subtly extending graphs of spaces} can be
applied inductively to produce a pre-cover $\hat{X}$ in which
every $\bar{f}_j$ extends to a full elevation.
\end{rem}

In the case where the vertex spaces of $X$ are simply connected,
note that the construction of proposition \ref{Subtly extending
graphs of spaces} preserves finiteness.

\begin{cor}\label{Subtly extending graphs}
Consider a graph of spaces $X$ with simply connected vertex spaces
and a combinatorial map $f:S^1\to X$. Suppose $\bar{X}\to X$ is a
finite-sheeted pre-covering.  Let $\bar{f}:\bar{S}^1\to \bar{X}$
be a diverse elevation of $f$ to $\bar{X}$.  Let
$$
\hat{S}^1\to S^1
$$
be a finite-sheeted covering of $S^1$ and suppose there exists an embedding $\bar{S}^1\hookrightarrow\hat{S}^1$ so that
$\hat{S}^1\to S^1$ extends $\bar{S}^1\to S^1$. Then there exists a
finite-sheeted pre-covering $\hat{X}\to X$ extending $\bar{X}\to
X$ so that $\bar{f}:\bar{S}^1\to\bar{X}$ extends to a full
elevation
$$
\hat{f}:\hat{S}^1\to\hat{X}
$$
of $f$.
\end{cor}

\subsection{The polygon decomposition}

It can be useful to treat compact surfaces with non-empty boundary
like graphs.  Given any such surface $\Sigma$ there exists a finite
collection of disjoint embedded arcs
$\alpha_1,\ldots,\alpha_n\subset\Sigma$, with
$\partial\alpha_i\subset\partial\Sigma$, so that the surface $V$
obtained by cutting along the $\alpha_i$ is a polygon with $4n$
sides. To see this, observe that for any such $\alpha_i$,
$$
\chi(\Sigma\smallsetminus\alpha_i)=\chi(\Sigma)+1.
$$
Taking the $\alpha_i$ as edge spaces and $V$ as a vertex space with
the obvious inclusions as edge maps decomposes $\Sigma$ as a graph
of spaces.  Scott used this decomposition in \cite{Sc1} to show that
the fundamental groups of such surfaces (that is, free groups) are
subgroup separable. We will call it the \emph{polygon decomposition}
of $\Sigma$, and will use it often.

\subsection{Ensuring diversity}

In the light of proposition \ref{Subtly extending graphs of spaces}
it will be useful to be able to impose diversity on collections of
elevations.

\begin{lem}\label{Imposing diversity} Let $\Sigma$ be a compact
surface with non-empty boundary and assume $\chi(\Sigma)<0$.
 Endow $\Sigma$ with the polygon decomposition described above.  Realize the boundary components as embeddings of circles $\{\delta_i:c_i\to\Sigma\}$.
 Let $\bar{\Sigma}\to\Sigma$ be a finite-sheeted pre-covering and fix a finite, disjoint collection $\{\bar{\delta}_j:\bar{c}_j\to\Sigma_k\}$ of
elevations of the $\delta_i$.  Then $\bar{\Sigma}\to\Sigma$
extends to a finite-sheeted pre-covering $\hat{\Sigma}\to\Sigma$
so that:
\begin{enumerate}
\item each component of $\hat{\Sigma}$ is a deformation retract of its pre-image in $\bar{\Sigma}$;
\item each elevation $\bar{\delta}_j$ extends to an elevation $\hat{\delta}_j:\hat{c}_j\to\hat{\Sigma}$;
\item the set $\{\hat{\delta}_j\}$ is diverse.  Indeed, whenever $x\in\partial\hat{c}_j$
and $y\in\partial\hat{c}_k$ are distinct, $\hat{\delta}_j(x)$ and
$\hat{\delta}_k(y)$ lie in distinct vertex spaces of $\hat{\Sigma}$.
\end{enumerate}
\end{lem}
\begin{pf}
Without loss, assume $\bar{\Sigma}$ is connected.  Let
$\tilde{\Sigma}$ be the canonical completion of $\bar{\Sigma}$.
Each elevation $\bar{\delta}_j$ extends canonically to some
coproduct of elevations
$\tilde{\delta}_j:\tilde{c}_j\to\tilde{\Sigma}$. Consider any pair
of connected components $\tilde{b}_j\subset\tilde{c}_j$
and $\tilde{b}_k\subset\tilde{c}_k$.  We will show that only finitely
many vertex spaces of $\tilde{\Sigma}$ intersect both
$\tilde{\delta}_j(\tilde{b}_j)$ and $\tilde{\delta}_k(\tilde{b}_k)$.

Suppose there are infinitely many vertex spaces of $\tilde{\Sigma}$
that $\tilde{\delta}_j(\tilde{b}_j)$ and
$\tilde{\delta}_k(\tilde{b}_k)$ both intersect.  Let
$\phi:\tilde{\Sigma}\to\Gamma(\tilde{\Sigma})$ be the natural map to
the underlying graph of $\tilde{\Sigma}$.  Then
$\phi\circ\tilde{\delta}_j$ and $\phi\circ\tilde{\delta}_k$ are
proper maps $\R\to\Gamma(\tilde{\Sigma})$ that share infinitely many
vertices of $\Gamma(\tilde{\Sigma})$ in their images.  Since
$\pi_1(\Gamma(\tilde{\Sigma}))$ is finitely generated,
$\phi\circ\tilde{\delta}_j(b_j)\cap \phi\circ\tilde{\delta}_k(b_k)$
contains the image of some injective ray
$[0,\infty)\to\Gamma(\tilde{\Sigma})$. Since $\chi(\Sigma)<0$ the
vertex space $V$ of the polygon decomposition of $\Sigma$ is a
$4n$-sided polygon for $n>1$. Therefore no disjoint pair of
elevations of the $\delta_i$ that enter by the same edge space of
$V$ leave by the same edge space of $V$, so it follows that
$\tilde{\delta}_j(b_j)$ and $\tilde{\delta}_k(b_k)$ intersect
non-trivially.  But then $\tilde{\delta}_j|_{b_j}$ and
$\tilde{\delta}_k|_{b_k}$ are isomorphic elevations by lemma
\ref{Disjoint images for elevations of embeddings}.

Therefore, by expanding $\Gamma(\bar{\Sigma})$ to some finite
graph
$\Gamma(\bar{\Sigma})\subset\hat{\Gamma}\subset\Gamma(\tilde{\Sigma})$
and setting $\hat{\Sigma}=p^{-1}(\hat{\Gamma})$ the result
follows.
\end{pf}

\section{From pre-covers to covers}
\label{From pre-covers to covers}

\subsection{Stallings' Principle}

Our strategy for proving subgroup separability is to replace a
pre-cover by a finite-sheeted pre-cover, and then complete that
pre-cover to a finite-sheeted cover.  The first step is described in
section \ref{Section 4}.  For the second step, we use an idea
developed by John Stallings in the context of graphs, which we
therefore name in his honour.  Similar ideas were also used by
Hempel \cite{He76} and Scott \cite{Sc1}.

\begin{prop}[Stallings' Principle]\label{Stallings}
Let $\bar{X}\to X$ be a pre-covering with the following property:
for each edge space $e$ of $X$ with edge maps $\partial_\pm:e\to
V_\pm$, for each conjugacy class $\mathcal{D}$ of subgroups of
$\pi_1(e)$, there exists a bijection between the set of elevations
of $\partial_+$ to $\bar{X}$ of degree
$\mathcal{D}$ and the set of elevations of $\partial_-$ to $\bar{X}$ of degree $\mathcal{D}$.  Then $\bar{X}\to
X$ can be extended to a covering $\hat{X}\to X$, with the same set
of vertex spaces as $\bar{X}$.
\end{prop}
\begin{pf}
For each edge space $e$ of $X$ and for each degree $\mathcal{D}$
there is a bijection between the set of elevations of $\partial_+$
of degree $\mathcal{D}$ and the set of elevations of $\partial_-$
of degree $\mathcal{D}$.  So there is a bijection between the set
of hanging elevations of $\partial_+$ of degree $\mathcal{D}$ and
the set of hanging elevations of $\partial_-$ of degree
$\mathcal{D}$.  Fix such a bijection.  For each hanging elevation
$\bar{\partial}_+:\bar{e}\to\bar{V}_+$ of $\partial_+$ take
$\bar{e}$ as an edge space for $\hat{X}$ and for edge maps take
the hanging elevation $\bar{\partial}_+$ and the corresponding
hanging elevation $\bar{\partial}_-$.  The resulting pre-cover,
with these additional edge spaces, has no hanging elevations and
so is a cover.
\end{pf}

\begin{thm}\label{Extending pre-covers of graphs}
Let $X$ be any graph of spaces with simply connected vertices (for
example, a graph).  Then any finite-sheeted pre-cover of $X$ can
be extended to a finite-sheeted cover.
\end{thm}
\begin{pf}
Let $\bar{X}$ be a finite-sheeted pre-cover of $X$.  Let $N$ be the maximum
number of pre-images in $\bar{X}$ of a vertex space in $X$.  Create a new
pre-cover from $\bar{X}$ by adding disconnected vertex spaces so
that every vertex space of $X$ has $N$ pre-images.  Now every edge
map of $X$ has $N$ elevations to the pre-cover, all of the same
(trivial) degree; so, by proposition \ref{Stallings}, the
pre-cover can be extended to a cover.
\end{pf}

\begin{cor}[M. Hall \cite{Hall49}]\label{Hall's theorem}
Finitely generated free groups are subgroup separable.
\end{cor}
\begin{pf}
Realize a free group $F$ as the fundamental group of a graph $X$.
Consider a finitely generated subgroup $H\subset F$, and the
corresponding covering $X^H\to X$.  Let $\Delta\subset X^H$ be a
finite subgraph.  Expanding $\Delta$ if necessary, it can be assumed
that $\Delta$ is connected and carries $\pi_1(X^H)$.  But $\Delta$
is a pre-cover of $X$, so can be extended to a cover $\hat{X}$.
\end{pf}

Indeed, a similar argument gives the stronger theorem originally proved by
Hall.

\begin{cor}[M.~Hall \cite{Hall49}]\label{Strong Hall's Theorem}
If $F$ is a finitely generated free group and $H\subset F$ is a
finitely generated subgroup then there exists a finite-index
subgroup $F'\subset F$ containing $H$ so that
$$
F'=H*F''
$$
\end{cor}
\begin{pf}
Again, realize $F$ as the fundamental group of a graph $X$ and let $X^H\to
X$ be the covering corresponding to $H$.  Since $H$ is finitely generated
there exists a finite, connected subgraph $X'\subset X^H$ so that $\pi_1(X')=H$.
 But $X'\to X$ is a pre-covering, so can be extended to a genuine finite-sheeted
 covering $\hat{X}$.  Since $X'$ is a subgraph of $\hat{X}$, $H$ is a free
 factor in $F'=\pi_1(\hat{X})$.
\end{pf}

Stallings' Principle can also be used to prove that free products
of subgroup separable groups are subgroup separable.

\begin{thm}[Burns \cite{Bu} and Romanovskii \cite{Ro}]\label{Free products}
If $G_+$ and $G_-$ are subgroup separable then $G=G_+*G_-$ is
subgroup separable.
\end{thm}
\begin{pf}
Realize $G_\pm=\pi_1(V_\pm,v_\pm)$ and let $X$ be the quotient of
$$
V_-\sqcup [-1,+1]\sqcup V_+
$$
obtained by identifying $\pm1$ with $v_\pm$, so $G=\pi_1(X)$.  For
$H\subset G$ a finitely generated subgroup, let $X^H\to X$ be the
corresponding covering and consider a finite subcomplex
$\Delta\subset X^H$.  Since $H$ is finitely generated there exists
a sub-graph of spaces $X'\subset X^H$ with finite underlying graph
that carries $H$; enlarging $X'$ if necessary, it can be assumed
that $\Delta\subset X'$.

For a vertex space $V'$ of $X'$ covering $V_\pm$ in $X$, let
$\Delta_{V'}\subset V'$ be a finite subcomplex containing:
\begin{enumerate}
\item the images of any edge spaces of $X'$ adjoining $V'$;
\item $\Delta\cap V'$.
\end{enumerate}
Since $G_\pm$ are subgroup separable and $\pi_1(V')$ is finitely
generated there exists an intermediate, finite-sheeted covering
$$
V'\to\bar{V}\to V
$$
such that $\Delta_{V'}$ embeds in $\bar{V}$.  Replace each vertex
space $V'$ of $X'$ by the corresponding $\bar{V}$; the edge maps
$\partial'_\pm:e\to V'$ descend to maps $\bar{\partial}_\pm:e\to
\bar{V}$.  These new vertex spaces and edge maps give an
intermediate, finite-sheeted, connected pre-covering
$$
X'\to\bar{X}\to X
$$
into which $\Delta$ embeds.

Each vertex space of $\bar{X}$ covers either $V_+$ or $V_-$.  Let
$$
P_\pm=\sum_{\bar{V}\to V_{\pm}}\deg(\bar{V}\to V_\pm).
$$
Without loss of generality, assume $P_+\geq P_-$.  Consider the
pre-cover given by the disjoint union of $\bar{X}$ and $P_+-P_-$
copies of $V_-$.  Any edge space $e$ of $X$ is a point; the edge
maps are inclusions $\partial_\pm:e\to V_\pm$ and their elevations
have only one possible (trivial) degree.  Now $\partial_+$ has $P_+$
elevations to vertex spaces of the pre-cover, and $\partial_-$ has
$P_-+(P_+-P_-)$ elevations.  So, by proposition \ref{Stallings},
the pre-cover $\bar{X}$ can be extended to a finite-sheeted cover $\hat{X}$.
It is connected, so the result follows.
\end{pf}

Stallings' Principle gives remarkable flexibility in constructing
covers of graphs and surfaces with boundary.  Here is an example
that will prove useful.  Note that the degree of an elevation of a
map $S^1\to X$ is (a conjugacy class of) a subgroup of $\Z$, and so
in particular corresponds uniquely to a non-negative integer.

\begin{cor}\label{Useful corollary}
Let $\Sigma$ be a compact surface with non-empty boundary.  Realize the boundary
components as embeddings of circles
$$
\delta_i: c_i\to\Sigma.
$$
Let $H\subset\pi_1(\Sigma)$ be a finitely generated subgroup, let
$\Sigma^H\to\Sigma$ be the corresponding covering and let
$\Delta\subset\Sigma^H$ be a compact subcomplex.  Let $\delta^H_j:
c^H_j\to\Sigma^H$ be all the finite-degree elevations of the $\delta_i$
 to $\Sigma^H$. Fix a finite collection of infinite-degree elevations
$$
\epsilon^H_k:c^H_k\to\Sigma^H
$$
of the $\delta_i$.  Then for all sufficiently large positive
integers $d$ there exists an intermediate finite-sheeted covering
$$
\Sigma^H\to\hat{\Sigma}\to\Sigma
$$
so that:
\begin{enumerate}
\item $\Delta\subset\hat{\Sigma}$;
\item the $\delta^H_j$ and $\epsilon^H_k$ all descend to distinct elevations
to $\hat{\Sigma}$;
\item if $\delta^H_j$ descends to $\bar{\delta}_j$ then $\deg(\bar{\delta}_j)=\deg(\delta^H_j)$;
\item if $\epsilon^H_k$ descends to $\bar{\epsilon}_k$ then $\deg(\bar{\epsilon}_k)=d$.
\end{enumerate}
\end{cor}
\begin{pf}
Consider the polygon decomposition of $\Sigma$, with vertex space
the polygon $V$ and arcs $\{e\}$ for edge spaces. The cover
$\Sigma^H$ inherits a graph-of-spaces decomposition. Since $H$ is
finitely generated, there exists a sub-graph of spaces
$\Sigma'\subset\Sigma^H$ with finite underlying graph such that the
inclusion $\Sigma'\hookrightarrow\Sigma^H$ is $\pi_1$-surjective.
Enlarging $\Sigma'$ if necessary, it may be assumed that:
\begin{enumerate}
\item $\Sigma'$ contains $\Delta$;
\item $\Sigma'$ contains the images of all the $\delta^H_j$;
\item for each $k$, the intersection $\mathrm{im}\epsilon^H_k\cap\Sigma'$
is a (non-empty) arc;
\item the set $\{\epsilon'_k\}$ is diverse.  (We can ensure this by lemma \ref{Imposing diversity}.)
\end{enumerate}
Note that $\Sigma'$ is a pre-cover of $\Sigma$ and $\Sigma^H$ is
its canonical completion. An elevation
$$
\epsilon^H_k:c^H_k\to \Sigma^H
$$
restricts to an elevation
$$
\epsilon'_k:c'_k\to\Sigma',
$$
where $c'_k$ is a compact subarc of $c^H_k\cong\R$.  If $d$ is
sufficiently large then it is clear that the $d$-fold covering
$\bar{c}_k\to c_i$ extends $c'_k\to c_i$. Since $\{\epsilon'_k\}$
is diverse we can apply corollary \ref{Subtly extending graphs}
repeatedly to extend $\Sigma'\to\Sigma$ to a pre-covering
$\bar{\Sigma}\to\Sigma$ so that every $\epsilon'_k$ extends to a
full elevation $\bar{\epsilon}_k$ of degree $d$. Finally, extend
$\bar{\Sigma}$ to a genuine cover $\hat{\Sigma}$ by theorem
\ref{Extending pre-covers of graphs}.
\end{pf}

\subsection{Homological assumptions}

In this subsection we use assumptions on the homology of the edge
spaces to make it possible to apply Stallings' Principle.

\begin{prop}\label{One-vertex pre-covers to covers}
Let $X$ be a graph of spaces with one vertex space $V$.  Suppose
every edge space $e$ is a circle and that, furthermore,
$H_1(e,\Z)$ is an infinite direct factor in $H_1(X,\Z)$.  Let
$\bar{X}\to X$ be a finite-sheeted pre-covering. Then there exist
finite-sheeted coverings $\hat{X}\to\bar{X}$ and $\tilde{X}\to X$ and an inclusion $\hat{X}\hookrightarrow\tilde{X}$ so that
\[ \begin{CD}
{\hat{X}} @>>>  {\tilde{X}} \\
@VVV    @VVV \\
 {\bar{X}} @>>>  X\\
\end{CD}\]
commutes.
\end{prop}
\begin{pf}
For each elevation $\bar{\partial}_\pm:\bar{e}\to\bar{V}$ of an
edge map $\partial_\pm:e\to V$ of $X$ to a vertex space $\bar{V}$
of $\bar{X}$, the degree $\deg({\bar{e}})$ of the
elevation can be thought of as a positive integer.  Let
$$
M=\prod_{\bar{e}}\deg{\bar{e}}
$$
where the product ranges over all elevations to $\bar{X}$ of edge maps of $X$.  Consider the finite-sheeted covering $X_M\to X$ corresponding to the kernel of the composition
$$
\pi_1(X)\to H_1(X)\to H_1(X,\Z/M\Z).
$$
Pull this covering back along $\bar{X}\to X$ to give a finite-sheeted
covering
$$
\hat{X}\to \bar{X}
$$
with the property that every elevation of an edge map of $X$ to a
vertex space of $\hat{X}$ is of degree $M$.

Consider
$$
P=\sum_{\hat{V}}\deg(\hat{V}\to V)
$$
where the sum ranges over all the vertex spaces of $\hat{X}$.  Any
edge map $\partial_\pm:e\to V$ has $P/M$ elevations to $\hat{X}$.  Since this is a constant, we can apply Stallings' Principle to deduce the result.
\end{pf}

This conclusion will be enough to deduce separability.  As an
application, we give a purely topological proof of Scott's theorem
that the fundamental groups of closed surfaces are subgroup
separable (cf.~\cite{Gi1} and \cite{Wi}).

\begin{cor}[Scott \cite{Sc1,Sc2}]\label{Surface groups}
Surface groups are subgroup separable.
\end{cor}
\begin{pf}
Let $\Sigma$ be a closed surface of negative Euler characteristic.  (Otherwise
the result is easy.) Therefore $\Sigma$ has a non-separating simple closed curve $\gamma$.  Note that $\gamma$ is primitive in $H_1(\Sigma)$.  Let $\Sigma_0$ be the compact surface with boundary that results
from deleting a small cylindrical neighbourhood of $\gamma$.  This
realizes $\Sigma$ as a graph of spaces, with one vertex space
$\Sigma_0$ and one edge space $e\cong S^1$. Consider any finitely
generated subgroup $H\subset\pi_1(\Sigma)$, and $\Sigma^H\to
\Sigma$ the corresponding covering. Fix a representative of a
curve $g\notin H$; the lift of $g$ to $\Sigma^H$ is not closed.

The cover $\Sigma^H$ inherits a graph-of-spaces decomposition from
$\Sigma$.  Since $H$ is finitely generated, there exists
some sub-graph of spaces $\Sigma'\subset\Sigma^H$ with finite
underlying graph that carries the whole of $H$.  Enlarging
$\Sigma'$ if necessary, it can be assumed that $\Sigma'$ contains
the lift of $g$. Note that $\Sigma'$ is a pre-cover of $\Sigma$.
Our first aim is to replace it by a finite-sheeted pre-cover.

Fix a large positive integer $d$.  By corollary \ref{Useful
corollary}, every vertex space $\Sigma'_0$ of $\Sigma'$ has an
intermediate covering
$$
\Sigma'_0\to\bar{\Sigma}_0\to\Sigma_0
$$
in which all the compact edge spaces have unchanged, finite
degree, all the non-compact edge spaces have degree $d$ and,
furthermore, the components of $g$ that cross $\Sigma'_0$ embed in
$\bar{\Sigma}_0$.  Replacing each $\Sigma'_0$ by $\bar{\Sigma}_0$
and assigning the edge maps to be the elevations that descend to
the new vertex groups gives a finite-sheeted pre-covering
$\bar{\Sigma}\to \Sigma$ containing a lift of $g$ that isn't
closed. Applying proposition \ref{One-vertex pre-covers to covers}
to $\bar{\Sigma}$ gives a finite-sheeted covering
$$
\hat{\Sigma}\to\bar{\Sigma}
$$
that extends to a finite-sheeted covering
$$
\tilde{\Sigma}\to\Sigma.
$$
It follows that the end-point of the lift of $g$ to
$\tilde{\Sigma}$ doesn't coincide with the end-point of the lift
of any $h\in H$; otherwise, $g$ would lift to a closed loop in
$\bar{\Sigma}$.
\end{pf}

We will also want to consider graphs of groups with two vertex
spaces.

\begin{prop}\label{Two-vertex pre-covers to covers}
Let $X$ be a graph of spaces with two vertex spaces $U$ and $V$.
Suppose every edge space $e$ is a circle and that, furthermore,
$H_1(e,\Z)$ is an infinite direct factor in $H_1(X,\Z)$.  Let
$\bar{X}\to X$ be a finite-sheeted, connected pre-covering. Then there exist
finite-sheeted coverings $\hat{X}\to\bar{X}$ and $\tilde{X}\to X$ and an inclusion $\hat{X}\hookrightarrow\tilde{X}$ so that
\[ \begin{CD}
{\hat{X}} @>>>  {\tilde{X}} \\
@VVV    @VVV \\
 {\bar{X}} @>>>  X\\
\end{CD}\]
commutes.
\end{prop}
\begin{pf}
As in the proof of proposition \ref{One-vertex pre-covers to
covers}, the degree $\deg(\bar{e})$ of an elevation
$\bar{\partial}_\pm:\bar{e}\to \bar{U}$ of an edge map
$\partial_\pm:e\to U$ of $X$ can be thought of as a positive
integer. Let
$$
M=\prod_{\bar{e}}\deg\bar{e}
$$
where the product ranges over all elevations of edge maps of $X$
to vertex spaces of $\bar{X}$. Define $X_M\to X$ to be the
covering associated to the kernel of the homomorphism
$$
\pi_1(X)\to H_1(X,\Z/M\Z).
$$
Set $\hat{X}\to\bar{X}$ to be the pullback of $X_M\to X$ along $\bar{X}\to
X$.  The cover $\hat{X}$ inherits a graph-of-spaces decomposition from $\bar{X}$.
Any edge space $\hat{e}\to e$ of $\hat{X}$ has degree $M$.

Let $U$ and $V$ be the two vertex spaces of $X$. Write
$$
P=\sum_{\hat{U}\to U}\deg(\hat{U}\to U)
$$
and
$$
Q=\sum_{\hat{V}\to V}\deg(\hat{V}\to V),
$$
where the sums range over all vertex spaces of $\hat{X}$ covering $U$ and
$V$ respectively.  Without loss of generality, $P\geq Q$.  Let
$$
\tilde{V}\to V
$$
be the covering associated to the kernel of the homomorphism
$$
\pi_1(V)\to H_1(V,\Z/M\Z).
$$
Let $N=|H_1(V,\Z/M\Z)|=\deg(\tilde{V}\to V)$.

Consider the pre-cover of $X$ consisting of the disjoint union of
$P-Q$ copies of $\tilde{V}$ and $N$ copies of $\hat{X}$.  Any edge
map $\partial_\pm:e\to U\sqcup V$ of $X$ has $NP/M$ elevations to
this pre-cover.  Since this is a constant, the pre-cover can be
extended to a genuine cover $\tilde{X}$.
\end{pf}

\section{Making pre-covers finite-sheeted}\label{Section 4}

As explained at the beginning of section \ref{From pre-covers to covers}, our strategy is to replace pre-covers by
finite-sheeted pre-covers and then complete these to genuine
finite-sheeted covers.  This section is devoted to a theorem that
shows, in certain circumstances, how to replace pre-covers by
finite-sheeted pre-covers.

\begin{thm}\label{Making pre-covers finite}
Let $X$ be a graph of spaces with two vertex spaces, constructed
as follows.
\begin{enumerate}
\item One vertex space, $L$, has subgroup separable fundamental group.
\item The other vertex space, $\Sigma$, is a compact surface with
non-empty boundary and $\chi(\Sigma)<0$.
\item Each edge space $e$ is a circle.  Of the edge maps:
$$
\partial_+:e\to L
$$
is an embedding; while
$$
\partial_-:e\to\Sigma
$$
is required to identify $e$ homeomorphically with a boundary component of
$\Sigma$.  Each boundary component of $\Sigma$ is identified with at most
one edge space of $X$.
\end{enumerate}
If $X'\to X$ is a pre-covering with finite underlying graph and
$\Delta\subset X'$ is a finite subcomplex then there exists an
intermediate finite-sheeted pre-covering
$$
X'\to\bar{X}\to X
$$
into which $\Delta$ embeds.
\end{thm}
\begin{pf}
Consider the polygon decomposition of $\Sigma$, so a covering space $\Sigma'\to\Sigma$ inherits a graph-of-spaces decomposition.

For each vertex $\Sigma'$ of $X'$ that lies above $\Sigma$, let
$\Delta_{\Sigma'}\subset\Sigma'$ be a sub-graph of spaces with finite, connected
underlying graph that carries $\pi_1(\Sigma')$ so that, furthermore:
\begin{enumerate}
\item $\Delta\cap\Sigma'\subset\Delta_{\Sigma'}$;
\item if $e'$ is an incident edge space of $X'$ then $\partial'_-(e')\cap\Delta_{\Sigma'}$
is non-empty;
\item the collection $\{\partial'_-\}$ of incident edge maps restricts to
a diverse collection of elevations to $\Delta_{\Sigma'}$.  (This
can be ensured by lemma \ref{Imposing diversity}.)
\end{enumerate}

For each vertex space $L'$ of $X'$ covering $L$ fix a finite, connected
subcomplex $\Delta_{L'}\subset L'$ such that:
\begin{enumerate}
\item $\Delta\cap L'\subset\Delta_{L'}$;
\item $\partial'_+((\partial'_-)^{-1}(\Delta_{\Sigma'}))\subset\Delta_{L'}$
whenever
$$
L'\stackrel{\partial'_+}{\leftarrow} e'\stackrel{\partial'_-}{\rightarrow}\Sigma'
$$
is an incident edge space.
\end{enumerate}
Since the edge maps $\partial_{\pm}$ of $X$ are embeddings, distinct edge
maps of $X'$ have disjoint images by lemma \ref{Disjoint images for elevations of embeddings}.

Since $\pi_1(L)$ is subgroup separable, for each vertex $L'\to L$ of $X'$
there exists an intermediate, finite-sheeted covering
$$
L'\to\bar{L}\to L
$$
so that $\Delta_{L'}$ injects into $\bar{L}$.

To construct $\bar{X}$, we first construct the underlying graph
$\bar{\Gamma}$.  Let $\Gamma'$ be the underlying graph of $X'$. We
will define an equivalence relation on the edges of $\Gamma'$.
For each vertex space $L'$ of $\Gamma$ covering $L$, let $\bar{L}$
be the finite-sheeted cover of $L$ fixed in the previous
paragraph.  Consider a pair of edges $e'_1$ and $e'_2$ of $X'$ that cover the same edge $e$ of $X$.  Let $\partial_+:e\to L$ be an edge map
of $e$; suppose $e'_1$ and $e'_2$ both adjoin the same vertex space $L'$
of $X'$ with edge maps $\partial^i_+:e'_i\to L'$.  If the edge maps $\partial^i_+$
descend to the same elevation $\bar{\partial}_+:\bar{e}\to\bar{L}$ then we
write $e'_1\sim e'_2$.  Let $\bar{\Gamma}=\Gamma'/\sim$.  Each vertex
of $\bar{\Gamma}$ is labelled $L$ or $\Sigma$, depending on which vertex
space of $X$ its pre-images in $\Gamma'$ cover.

Each $L$-vertex $\bar{v}$ of $\bar{\Gamma}$ corresponds to a unique
$L$-vertex $v'$ of $\Gamma'$.  If the vertex space corresponding
to $v'$ is $L'$, then the vertex space of $\bar{v}$ is $\bar{L}$.

Each $\Sigma$-vertex $\bar{u}$ of $\bar{\Gamma}$ corresponds to an
equivalence class of finitely many vertices $u'_i$ of $\Gamma'$,
with corresponding vertex spaces $\Sigma'_i$.  Associate to
$\bar{u}$ the pre-covering
$$
\Sigma''=\sqcup_i \Delta_{\Sigma'_i}\to\Sigma.
$$
Later, we will complete $\Sigma''$ to a cover $\bar{\Sigma}$.

Each edge of $\bar{\Gamma}$ corresponds to a finite equivalence
class of edges $e'_i$ of $\Gamma'$, all covering the same edge $e$
of $X$.   The edge maps of $e$ are $\partial_+:e\to L$ and
$\partial_-:e\to\Sigma$, and the edge maps of $e'_i$ are
elevations $\partial_+^i:e'_i\to L'$ and $\partial_-^i:e'_i\to
\Sigma'_i$.  We aim to find suitable edge spaces $\bar{e}$ and
edge maps $\bar{\partial}_+:\bar{e}\to\bar{L}$ and
$\bar{\partial}_-:\bar{e}\to\bar{\Sigma}$.

By the definition of the equivalence relation on the edges of $\Gamma'$, every
$\partial^i_+$ descends to the same elevation
$$
\bar{\partial}_+:\bar{e}\to\bar{L}.
$$
This gives the edge space and one of the edge maps.

On the other side, consider
$$
e''=\bigsqcup_i (\partial^i_-)^{-1}(\Delta_{\Sigma'_i})\subset
\sqcup_i e'_i.
$$
The restriction
$$
\partial''_-=\sqcup_i \partial^i_- |_{e''}
$$
is an elevation of $\partial_-$ to $\Sigma''$.

There is a natural map $\iota:e''\to\bar{e}$, namely the coproduct of the
restrictions of the covering maps $e'_i\to\bar{e}$.  The claim is that $\iota$ is an injection.  We can then apply corollary \ref{Subtly extending graphs}
to extend $\Sigma''$ in such a way that $\partial''_-$ extends to $\bar{e}$.

Consider $t_1,t_2\in e''$ and suppose that $\iota(t_1)=\iota(t_2)$.  Assume
that $e''_j\subset e''$ is the connected component containing $t_j$. Since $e''=\sqcup_i(\partial^i_-)^{-1}(\Delta_{\Sigma'_i})$ it follows that $\partial^j_+(t_j)\in\Delta_{L'}$.
Since $L'\to\bar{L}$ is injective on $\Delta_{L'}$ it follows that $\partial^1_+(t_1)=\partial^2_+(t_2)$.
 But then $\partial^1_+$ and $\partial^2_+$ are elevations of $\partial_+$
 with non-disjoint image, so they are in fact isomorphic and $e'_1=e'_2$.
 Since $\partial^1_+:e'_1\to L'$ is an elevation of an injective map it is
 itself injective, so it now follows that $t_1=t_2$.  This proves the claim.

To complete the proof of the theorem, we now apply corollary
\ref{Subtly extending graphs} (invoking remark \ref{Diversity is
useful}) to extend every $\Sigma''$ so that the elevations
$\partial''_-:e''\to \Sigma''$ extend to full elevations
$\bar{\partial}_-:\bar{e}\to\Sigma''$.  Finally, by theorem
\ref{Extending pre-covers of graphs}, each pre-cover $\Sigma''$ can
be extended to a genuine finite-sheeted cover $\bar{\Sigma}$.
\end{pf}

\begin{cor}\label{Main technical result}
Let $X$ be a graph of spaces satisfying the hypotheses of theorem
\ref{Making pre-covers finite} and such that, furthermore, for
each edge space $e$ of $X$, $H_1(e,\Z)$ is an infinite direct
factor of $H_1(X,\Z)$.  Then $\pi_1(X)$ is subgroup separable.
\end{cor}
\begin{pf}
Fix a base-point $x\in X$.   Let $H\subset\pi_1(X)$ be a finitely
generated subgroup and fix a representative of
$g\in\pi_1(X)\smallsetminus H$. Let $(X^H,x')\to (X,x)$ be the
covering corresponding to $H$. Since $H$ is finitely generated,
there exists a sub-graph of spaces $X'\subset X^H$ with finite
underlying graph so that $X'$ contains the (based) lift $g'$ of the
representative of $g$.  By theorem \ref{Making pre-covers finite}
there exists a finite-sheeted, intermediate pre-cover
$$
(X',x')\to(\bar{X},\bar{x})\to (X,x)
$$
into which the image of $g'$ embeds.  By proposition
\ref{Two-vertex pre-covers to covers}, there exists a
finite-sheeted covering
$$
\hat{X}\to \bar{X}
$$
that extends to a finite-sheeted covering
$$
\tilde{X}\to X.
$$
As in the proof of corollary \ref{Surface groups}, it follows that the end-point
of the lift of the representative of $g$ to $\tilde{X}$ doesn't coincide
with the end-point of the lift of any representative of any $h\in H$.
\end{pf}

\section{Elementarily free groups}

\subsection{$\omega$-residually free towers}

Elementarily free groups are an important class of limit groups.
For an introduction to the theory of limit groups, see
\cite{BF03}.  The simplest definition of a limit group is in terms
of the property of being $\omega$-residually free.

\begin{defn}
Fix a finitely generated free group $F$ of rank at least 2.  A
group $G$ is \emph{residually free} if, for any $g\in G$, there
exists a homomorphism $f:G\to F$ so that $f(g)\neq 1$.  A finitely
generated group $G$ is a \emph{limit group}, or
\emph{$\omega$-residually free}, if for any finite subset
$S\subset G$ there exists a homomorphism $f:G\to F$ such that
$f|_S$ is injective.
\end{defn}

Attempts to solve Tarski's problems on the elementary theory of free groups
have led to a comprehensive structure theory for limit groups, most easily stated in terms of $\omega$-residually free towers.

\begin{defn}\label{Towers}
A \emph{tower space of height $0$}, denoted $X_0$, is a one-point
union of finitely many compact graphs, tori, and closed hyperbolic
surfaces of Euler characteristic less than -1.

A \emph{tower space of height $h$}, denoted $X_h$, is obtained
from a tower space $X_{h-1}$ of height $h-1$ by attaching one of
two sorts of blocks.
\begin{enumerate}
\item \textbf{Quadratic block.}  Let $\Sigma$ be a connected
compact hyperbolic surface with boundary, with each component
either a punctured torus or having $\chi\leq-2$.  Then $X_h$ is
the quotient of $X_{h-1}\sqcup\Sigma$ obtained by identifying the
boundary components of $\Sigma$ with curves on $X_{h-1}$, in such
a way that there exists a retraction $\rho:X_h\rightarrow
X_{h-1}$. The retraction is also required to satisfy the property
that $\rho_*(\pi_1(\Sigma))$ is non-abelian.

\item \textbf{Abelian block.} Let $T$ be an $n$-torus, and fix
a coordinate circle $c$.  Fix a loop $\gamma$ in $X_{h-1}$ that
generates a maximal abelian subgroup in $\pi_1(X_{h-1})$.  Then
$X_h$ is the quotient of $X_{h-1}\sqcup (S^1\times[0,1])\sqcup T$
obtained by identifying $S^1\times\{0\}$ with $\gamma$ and
$S^1\times\{1\}$ with $c$.
\end{enumerate}
A tower space is called \emph{hyperbolic} if no tori are used in
its construction.
\end{defn}

\begin{defn}
An \emph{($\omega$-residually free) tower of height $h$}, denoted
$L_h$, is the fundamental group of a tower space of height $h$.
\end{defn}

A tower space $X_h$ has a natural graph-of-spaces decomposition
$\Gamma_X$ with two vertex spaces, namely $X_{h-1}$ and the block
at height $h$; the edge spaces are circles.  By the Seifert--van
Kampen Theorem, towers have a corresponding graph-of-groups
decomposition.

The following deep theorem of Sela (see \cite{Se6}) will, for our
purposes, serve as a definition of elementarily free groups.

\begin{thm}
A group is \emph{elementarily free} if and only if it is the
fundamental group of a hyperbolic tower space.
\end{thm}

Towers are limit groups. Another theorem of Sela \cite{Se2} and,
independently, O.~Kharlampovich and A.~Myasnikov \cite{KM98b}, shows
how towers provide a structure theory for limit groups.

\begin{thm}\label{Structure of limit groups}
A group is a \emph{limit group} if and only if it is a finitely
generated subgroup of an $\omega$-residually free tower.
\end{thm}

One consequence of this theorem is that all abelian subgroups of limit groups are finitely generated.  For elementarily free groups, one can do better.

\begin{lem}[Cf. Lemma 2.1 of \cite{Se1}]\label{Abelian subgroups of EFGs}
Every abelian subgroup of an elementarily free group is cyclic.
\end{lem}

A key feature of the definition of a tower is the retraction
$\rho:X_h\rightarrow X_{h-1}$.  In the abelian case, the
retraction simply projects $T$ onto the coordinate circle $c$, and
thence to $\gamma$.  In both cases, $\rho$ induces a retraction
$\rho_*:L_h\rightarrow L_{h-1}$ on the level of fundamental
groups.

We will often use the retraction to pull finite covers back from
$X_{h-1}$ to $X_h$. It is worth noting that such pullbacks inherit
from $X_h$ a similar graph-of-spaces decomposition.

\begin{lem}\label{Pullback covers}
Let $X$ be a complex with a graph-of-spaces decomposition
$\Gamma_X$, such that there is a retraction $\rho:X\rightarrow V$
to a vertex space.  Let $\hat{V}\rightarrow V$ be a connected covering
of degree $d$, and let $\hat{X}\rightarrow X$ be the connected covering
obtained by pulling back along $\rho$.  Then:
\begin{enumerate}
\item $\hat{X}\rightarrow X$ is of degree $d$ and inherits a graph-of-spaces decomposition $\hat{\Gamma}_X$;
\item the pre-image of $V$ in $\hat{X}$ is a (connected) vertex
space of $\hat{\Gamma}_X$ homeomorphic to $\hat{V}$;
\item $\hat{X}\rightarrow X$ extends $\hat{V}\rightarrow V$ and $\hat{X}$
inherits a retraction to $\hat{V}$ covering $\rho$.
\end{enumerate}
\end{lem}

\subsection{Simplifying the gluing maps}

The aim is to apply corollary \ref{Main technical result} to prove
that elementarily free groups are subgroup separable.  We
therefore need to simplify the gluing maps of the tower spaces by
passing to finite-sheeted covers so that the hypotheses of
corollary \ref{Main technical result} are satisfied.  First, we
show how to satisfy the homological conditions.

\begin{lem}\label{Satisfying homological conditions}
If $X$ is a connected complex with $L=\pi_1(X)$ residually free and
$\gamma:S^1\to X$ is a loop in $X$ then there exists a
finite-sheeted covering $\bar{X}\rightarrow X$ so that every
elevation of $\gamma$ to $\bar{X}$ generates an infinite direct
factor in $H_1(\bar{X})$.
\end{lem}
\begin{pf}
Fix a base-point in $X$ and without loss of generality assume that
$\gamma$ is a based loop representing an element of $L$. Since $L$
is residually free, there exists a homomorphism $f:L\rightarrow F$
with $f(\gamma)\neq 1$. By corollary \ref{Strong Hall's Theorem},
there exists a finite-index subgroup $F'\subset F$ containing
$f(\gamma)$ such that $\langle f(\gamma)\rangle$ is a free factor in
$F'$; in particular, $\langle f(\gamma)\rangle$ is a direct factor
in $H_1(F')$. Let $\bar{X}\rightarrow X$ be the covering
corresponding to the subgroup $f^{-1}(F')$. Every elevation
$\bar{\gamma}$ of $\gamma$ to $\bar{X}$ corresponds to a conjugate
of a power of $\gamma\in L$. Since $f(\gamma)$ is primitive and has
infinite order in $H_1(F')$, it follows that $\bar{\gamma}$
generates an infinite direct factor in $H_1(\bar{X})$.
\end{pf}

It remains to ensure that the gluing maps are injective.  It will
prove useful that cyclic subgroups of towers are separable.  This
can easily be proved by induction on height, but it also follows
from the stronger result that maximal abelian subgroups of limit
groups are \emph{closed in the pro-free topology}.

\begin{lem}\label{Maximal abelian subgroups are separable}
Let $L$ be residually free, $A\subset L$ a maximal abelian subgroup
and $g\in L\smallsetminus A$.  Then there exists a homomorphism
$f:L\to F$ so that $f(g)\notin f(A)$.

Furthermore, if $L$ is an elementarily free group, $A\subset L$ is
any cyclic subgroup and $g\in L\smallsetminus A$ then there exists $f:L\to F$
with $f(g)\notin f(A)$.

It follows by theorem \ref{Hall's theorem} that $A$ is
separable.
\end{lem}
\begin{pf}
Choose some $a\in A$ such that $[g,a]\neq 1$.  Then there exists a
homomorphism $f:A\to F$ so that $f([g,a])\neq 1$.  In particular,
$f(g)\notin f(A)$.

If $L$ is elementarily free then by lemma \ref{Abelian subgroups of
EFGs} $A$ is contained in some maximal abelian, cyclic subgroup
$A'\subset L$ and we are reduced to the case $g\in A'\smallsetminus
A$.  Let $a'$ generate $A'$.  Then $f:L\to F$ with $f(a')\neq 1$ is
as required.
\end{pf}

We will use this property to desingularize the gluing curves, but
first we need a little geometry of non-positive curvature. For the
definition of and general references on \emph{CAT(0) geodesic metric
spaces} see \cite{BH99}.  Such spaces are contractible and geodesics
are unique (corollary II.1.5 and proposition II.1.4 of \cite{BH99}).
A metric space in which every point has a CAT(0) neighbourhood is
called \emph{non-positively curved}.  We will make use of the
following fact, which follows from the results of pages 229 to 231
of \cite{BH99}.

\begin{lem}\label{Non-positive curvature}
Let $\gamma:S^1\to X$ be a loop in a compact, non-positively curved
metric space $X$.  Then $\gamma$ is freely homotopic to a local
geodesic $\gamma_0$.  The length of $\gamma_0$ is minimal in the
homotopy class of $\gamma$.  It follows that, if $l(\gamma)$ denotes
the length of $\gamma$,
$$
l(\gamma_0^n)=|n|l(\gamma_0)
$$
for any integer $n$.
\end{lem}

That tower spaces can be endowed with non-positively curved metrics
follows immediately from, for example, proposition II.11.13 of
\cite{BH99}. More generally, Emina Alibegovic and Mladen Bestvina
\cite{AB04} showed that all limit groups act geometrically on a
CAT(0) space with the isolated flats property.

\begin{lem}\label{Curve lifting}
Let $X$ be a compact, connected, non-positively curved metric space
such that cyclic subgroups of $\pi_1(X)$ are separable and suppose
$\gamma:S^1\to X$ is a homotopically non-trivial continuous map.
Then, after modifying $\gamma$ by a homotopy, there exists a
finite-sheeted covering $\hat{X}\to X$ such that every elevation of
$\gamma$ to $\hat{X}$ is injective.
\end{lem}
\begin{pf}
Since $X$ is compact we can assume that balls of radius less than
some $\epsilon>0$ are CAT(0).  After a homotopy, $\gamma$ can be
assumed to be a local isometry by lemma \ref{Non-positive
curvature}.  Parametrizing by arc length we view $\gamma$ as a map
$[0,l]\to X$ with $\gamma(0)=\gamma(l)$. Let $D(\gamma)\subset
[0,l]$ consist of those $t\in[0,l]$ such that there exists $s<t$
with $\gamma(s)=\gamma(t)$. Since balls of radius $\epsilon>0$ are
CAT(0), if $\gamma(s)=\gamma(t)$ then either $s=t$ or
$|s-t|\geq\epsilon$.

Note that $D(\gamma)$ is closed.  For, suppose $t_n\in
D(\gamma)$ and $t_n\to t$.  For each $n$ there exists $s_n\leq t_n-\epsilon$
with $\gamma(s_n)=\gamma(t_n)$.  Passing to a subsequence, the $s_n$ converge
to some $s\leq t-\epsilon$ with $\gamma(s)=\gamma(t)$.

Set $t_0=\min D(\gamma)$ and without loss of generality take
$\gamma(t_0)$ as a base-point for $\pi_1(X)$.  Then
$\gamma=\gamma_1\gamma_2$ where $\gamma_1$ is an embedded closed
path.  The length of $\gamma_1$ is less than the length of $\gamma$,
which is the shortest element of $\{\gamma^n|n\in\Z\}$ by lemma
\ref{Non-positive curvature}, so
$\gamma_1\notin\langle\gamma\rangle$. Using that
$\langle\gamma\rangle$ is separable, therefore, there exists a
finite-sheeted cover $X'$ of $X$ to which $\gamma$ lifts but
$\gamma_1$ doesn't.  In particular, $\min D(\gamma')>t_0$.
Furthermore, since the ball $B(\gamma'(t_0),\epsilon)\subset X'$ is
CAT(0) and isometric to $B(\gamma(t_0),\epsilon)\subset X$ it
follows that $\min D(\gamma')>t_0+\epsilon$.  Repeating inductively
gives a finite-sheeted covering $X'\to X$ and a lift $\gamma':S^1\to
X'$ with $D(\gamma')=\{l\}$; that is, $\gamma'$ is an embedding.

Let $\hat{X}\to X$ be the covering defined by
$$
\pi_1(\hat{X})=\bigcap_{g\in\pi_1(X)}\pi_1(X')^g,
$$
a finite-sheeted, regular covering.  Then every elevation
$\hat{\gamma}$ to $\hat{X}$ of $\gamma$ descends to $\gamma'$ after
composition with a covering automorphism of $\hat{X}$.  Since
elevations of embeddings are injective, it follows that
$\hat{\gamma}$ is injective.
\end{pf}

\begin{prop}\label{Clean tower}
Let $X_h$ be a tower space, constructed as above by attaching a
quadratic block $\Sigma$ to a space $X_{h-1}$ of height $h-1$. Then,
after modifying the gluing maps by a homotopy, there exists a
connected covering $Y_h\rightarrow X_h$ with an inherited
graph-of-spaces decomposition $\Gamma_Y$, with one vertex space
given by a connected cover $Y_{h-1}\rightarrow X_{h-1}$ and the
remaining vertex spaces given by connected coverings
$\Sigma_i\rightarrow\Sigma$, so that each edge space of $Y_h$
generates an infinite direct factor in $H_1(Y_h)$ and all the gluing
maps are injective.
\end{prop}
\begin{pf}
Let $\gamma_i$ be the gluing curves in $X_{h-1}$.  By lemma
\ref{Satisfying homological conditions}, for each $i$ there exists
a finite-sheeted covering $\bar{X}_i\to X_{h-1}$ so that every
elevation of $\gamma_i$ to $\bar{X}_i$ is primitive and of
infinite order in homology. By lemma \ref{Curve lifting}, for each
$i$ there exists a finite-sheeted covering $\hat{X}_i\to X_{h-1}$
so that every elevation of $\gamma_i$ to $\hat{X}_i$ is injective.
Let $Y_{h-1}\to X_{h-1}$ be the finite-sheeted covering
corresponding to the finite-index subgroup
$$
\bigcap_i \pi_1(\bar{X}_i)\cap \bigcap_i\pi_1(\hat{X}_i)\subset
\pi_1(X_{h-1}).
$$
Then the covering $Y_h\to X_h$ obtained by pulling $Y_{h-1}\to
X_{h-1}$ back along the retraction is as required.
\end{pf}

We are, finally, in a position to prove the main theorem.

\begin{thm}\label{Main theorem}
Every elementarily free group is subgroup separable.
\end{thm}
\begin{pf}
The proof is by induction on height.  Let $X_h$ be a
tower space of height $h$.  If $h=0$ then the result follows by
corollaries \ref{Hall's theorem} and \ref{Surface groups} and
theorem \ref{Free products}.  For the general case, let $X_h$ be constructed
as in definition \ref{Towers} by gluing a quadratic block $\Sigma$ to a
space $X_{h-1}$ of height $h-1$.  By induction, assume
$\pi_1(X_{h-1})$ is subgroup separable.  By proposition \ref{Clean
tower} there exists a finite-sheeted cover $Y_h$ of $X_h$
constructed by gluing surfaces $\Sigma_i$ to a finite-sheeted
cover $Y_{h-1}$ of $X_{h-1}$, and such that all the gluing maps
are injective and boundary components generate infinite direct
factors in homology.  Applying corollary \ref{Main technical
result} to each $\Sigma_i$ in turn it follows that $\pi_1(Y_h)$,
and hence $\pi_1(X_h)$, is subgroup separable.
\end{pf}

\section{A less direct proof}\label{Less direct}

In this section we outline a less direct proof.  In essence, we take a more
extreme finite-sheeted covering of the tower space that forces the attached
quadratic blocks to be of positive genus and then apply Gitik's theorem.
The arguments used by Bridson, Tweedale and myself in \cite{BTW} amount to the following theorem.

\begin{thm}[Proof of theorem 2.6 of \cite{BTW}]\label{BTW's result}
Let $\mathcal{C}$ be a class of finitely generated groups with the
following properties.
\begin{enumerate}
\item All free groups and hyperbolic surface groups lie in $\mathcal{C}$.
\item $\mathcal{C}$ is closed under taking subgroups.
\item If $H\in\mathcal{C}$ and $H$ is a finite-index subgroup of
$G$ then $G\in\mathcal{C}$.
\item If $G_1,G_2\in\mathcal{C}$ then
$G_1*G_2\in\mathcal{C}$.
\item Let $\Sigma$ be a positive-genus, hyperbolic surface with one boundary
component and $G\in\mathcal{C}$.  Then
$$
G*_{\pi_1(\partial\Sigma)}\pi_1(\Sigma)\in\mathcal{C}.
$$
\end{enumerate}
Then every elementarily free group lies in $\mathcal{C}$.
\end{thm}

In \cite{Gi1}, Gitik proved the following.

\begin{thm}[Gitik \cite{Gi1}]\label{Gitik's theorem}
If $G$ is subgroup separable, $F$ is free and $f\in F$ has no proper
roots, then $G*_{\langle f\rangle}F$ is subgroup separable.
\end{thm}

\begin{altpfof}{theorem \ref{Main theorem}}
Let $\mathcal{C}$ be the class of subgroup separable groups.
Property 1 of theorem \ref{BTW's result} is satisfied by corollaries
\ref{Hall's theorem} and \ref{Surface groups}.  It is immediate from
the definition that $\mathcal{C}$ is closed under taking subgroups.
Lemma \ref{Elementary properties} gives property 3 and theorem
\ref{Free products} gives property 4.  Finally, property 5 follows
from theorem \ref{Gitik's theorem}.  So, by theorem \ref{BTW's
result}, the class $\mathcal{C}$ contains all elementarily free
groups.
\end{altpfof}

\subsection*{Acknowledgement}

I would like to thank my supervisor Martin Bridson for all his
encouragement and helpful comments.

\bibliographystyle{plain}
\bibliography{lerf}

\bigskip\bigskip\centerline{\textbf{Author's address}}

\smallskip\begin{center}\begin{tabular}{l}%
Department of Mathematics,\\%
Imperial College,\\%
London SW7 2AZ.\\%
\\%
{\texttt{henry.wilton@imperial.ac.uk}}
\end{tabular}\end{center}

\end{document}